*Research Article*

# A Note on the Gaps in the Support of Discretely Infinitely Divisible Laws

## Anthony G. Pakes[1] and S. Satheesh[2]

[1] *School of Mathematics & Statistics, University of Western Australia, 35 Stirling Highway, Crawley, WA 6009, Australia*
[2] *Department of Applied Sciences, Vidya Academy of Science and Technology, Thalakkottukara, Thrissur 680 501, India*

Correspondence should be addressed to S. Satheesh; ssatheesh1963@yahoo.co.in

We discuss the nature of gaps in the support of a discretely infinitely divisible distribution from the angle of compound Poisson laws/processes. The discussion is extended to infinitely divisible distributions on the nonnegative real line.

## 1. Introduction

We begin with the following definitions in which $N$ denotes the set of positive integers and $Z_+$ the set of nonnegative integers.

*Definition 1.* A $Z_+$-valued random variable (r.v) $X$ is infinitely divisible (ID) if for every $n \in N$ there exists a r.v $X_n$ such that

$$X = X_{n,1} + \cdots + X_{n,n}, \qquad (1)$$

where the $X_{n,i}$'s are independent copies of $X_n$.

*Definition 2.* A $Z_+$-valued r.v $X$ is discretely infinitely divisible (DID) if for every $n \in N$ there exists a $Z_+$-valued r.v $X_n$ such that (1) is true.

The distinction between the notions of infinite divisibility (ID) of discrete distributions in general and discrete infinite divisibility (DID) is often not made and this sometimes causes confusion. For example, the geometric distribution on $\{1, 2, \ldots\}$ is ID but not DID whereas that on $Z_+ = \{0, 1, 2, \ldots\}$ is DID and thus is ID. Blurring these concepts may have lead to Remark 9 in Bose et al. [1] asserting that if a $Z_+$-valued ID distribution assigns a positive probability to the integer 1, then its support cannot have any gaps. Satheesh [2] gave the following simple example to show that this is not always true.

*Example 3.* Consider the r.v $X$ with probability generating function (PGF) $P_1(s) = ps/(1 - (1 - p)s^k)$, $k > 1$ integer.

Obviously, $P\{X = 1\} > 0$, but its support has gaps as the atoms are separated by $k - 1$ integers.

In the sequel supp($X$) or supp($F$) denotes the support of $X$ or that of its distribution function (d.f.) $F$. To emphasize what now is obvious, if $X$ is ID and $\ell_X = \min\{j \in Z_+ : P\{X = j\} > 0\}$ is its left extremity and if $k \in N$, then $k(X - \ell_X)$ is DID and its support has gaps if $k \geq 2$. Also, in terms of Definition 1, $\ell_{X_n} = \ell_X/n$, and hence supp($X_n$) = supp($X$) only if $X$ is DID. To avoid such cases, unless otherwise stated, we assume that supp($X - \ell_X$) is aperiodic; that is, its greatest common divisor is unity. This natural restriction implies the following definition.

*Definition 4.* The support supp($X$) of the ID discrete r.v $X$ on $Z_+$ contains gaps if it is a proper subset of $\ell_X + Z_+$. (Note that supp($X$) is necessarily an infinite set.)

The distinction between Definitions 1 and 2 is made clear in Grandell [3, page 26; a result due to Kallenberg], Satheesh [2], and Steutel and van Harn [4, page 23]. Implications of the r.v $X$ being DID include the following.

(i) $\ell_X = 0$.

(ii) Let $F_n$ be the $n$th convolution root of $F$; that is, $F_n$ is the d.f.of $X_n$ in (1). Then supp($X$) = supp($X_n$) or in terms of their d.f's supp($F$) = supp($F_n$) for all $n \in N$.

Table 1

| $N(t) = n = n_1 + n_2$ | 0 | 1 | 2 | 3 | 4 | 5 | $\cdots$ |
|---|---|---|---|---|---|---|---|
| $X(t) = n_1 y_1 + n_2 y_2$ | 0 | 1, 3 | 2, 4, 6 | 3, 5, 7, 9 | 4, 6, 8, 10, 12 | 5, 7, 9, 11, 13, 15 | $\cdots$ |

(iii) Let $P$ and $P_n$ denote the PGF's of $X$ and $X_n$, respectively. Then $P(s) = [P_n(s)]^n$, $|s| \leq 1$, for every $n \in N$, that is, for all $n \in N$; the $n$th root of a DID PGF is again a PGF.

(iv) $X$ is compound Poisson; that is, $P(s) = \exp\{-\lambda[1 - Q(s)]\}$, where $Q$ is a PGF and the rate parameter $\lambda > 0$.

We conclude this review by noting that Satheesh [2] modified Remark 9 in Bose et al. [1] as follows.

**Theorem 5.** *Let $X$ be a $Z_+$-valued ID r.v with $P\{X = 0\} > 0$. Then* supp$(X)$ *cannot have any gaps if and only if $P\{X = 1\} > 0$.*

The purpose of this note is to complement the discussion on the gaps in the support of DID laws in Satheesh [2] from the angle of compound Poisson laws/processes. Some general considerations about the supports of positive ID distributions are discussed in Section 3.

## 2. Compound Poisson Processes and DID

If we think in terms of compound Poisson laws/processes, the picture is clearer. Let $\{N(t), t \geq 0\}$ be a Poisson process with rate parameter $\lambda > 0$ and let $\{Y_i, i \geq 1\}$ be a sequence of independent copies of a generic jump random variable $Y$. Then, $X(t) = \sum_{i=1}^{N(t)} Y_i$ describes a compound Poisson process, that is, a subordinator with zero drift and Lévy measure which assigns mass $\lambda P(Y = j)$ to $\{j\}$. In keeping with the usual convention, we assume that $P(Y = 0) = 0$.

Suppose that the support of the jump distribution (of $Y$) is $\{y_1, y_2, \ldots, y_m\}$ with all $y_i$'s positive integers and $m \leq \infty$. From (iv) above one may identify the correspondence that $P_t(s) = \exp\{-\lambda t[1 - Q(s)]\}$ is the PGF of $X(t)$, and we set $X = X(1)$. It is clear that supp$(X(t)) = \{n_1 y_1 + n_2 y_2 + \cdots : n_i \in Z_+\}$ is independent of $t$. The following examples illustrate possible behaviors of gaps in supp$(X)$.

*Example 6.* Consider the case when $m = 2$, $y_1 = 1$, and $y_2 = 3$. Here supp$(X(t)) = Z_+$; there are no gaps. Although this is obvious, it can be seen better by considering the semigroup elements generated by successively increasing values of the Poisson driving process. So consider the unordered pairwise partitions $(n_1, n_2)$ of $n = 0, 1, \ldots$, respectively. The possible values of $X(t)$ corresponding to the first few values of $n$ are shown in Table 1.

*Example 7.* Next consider the case when $m = 2$, $y_1 = 2$, and $y_2 = 3$. As before we now have the possible outcomes in Table 2.

Thus, supp$(X(t)) = \{0, 2, \ldots\}$ where, here and below, the dots denote all integers larger than the last one written. Clearly, the support omits $\{1\}$, a gap.

Table 2

| $N(t) = n = n_1 + n_2$ | 0 | 1 | 2 | 3 | 4 | $\cdots$ |
|---|---|---|---|---|---|---|
| $X(t) = n_1 y_1 + n_2 y_2$ | 0 | 2, 3 | 4, 5, 6 | 6, 7, 8, 9 | 8, 9, 10, 11, 12 | $\cdots$ |

*Example 8.* If we consider the case when $m = 2$, $y_1 = 3$ and $y_2 = 7$ then see Table 3.

Thus, supp$(X(t)) = \{0, 3, 6, 7, 9, 10, 12, \ldots\}$, which omits $\{1, 2, 4, 5, 8, 11\}$, so several gaps.

*Example 9.* If we consider the case when $m = 2$, $y_1 = 4$, and $y_2 = 9$, then (see Table 4).

Thus supp$(X(t)) = \{0, 4, 8, 9, 12, 13, 16, 17, 18, 20, 21, 22, 24, \ldots\}$ and there are many gaps.

It is important to notice that the widths of successive gaps vary haphazardly but they do tend to become narrower as they approach the infinite region of the support containing all larger integers; that is, supp$(X(t))$ contains a translate of $Z_+$. In the last example, the gap lengths are $(4, 3, 2)$ and it illustrates that the gap lengths decrease. In general, we have the following result.

**Theorem 10.**

(i) *The support of the compound Poisson process $X(t)$ is $\cup\{n_1 y_1 + \cdots + n_m y_m; n_1, \ldots, n_m \in Z_+\}$; the additive semigroup generated by the $y_i$'s. That is,* supp$(\{X(t)\}) = $ *s.g.(supp$(Y_i)$).*

(ii) *Given $2 \leq n \in N$ the support relation* supp$(F_n) = $ supp$(F)$ *holds if and only if $X$ is DID, and then it holds for all $n$.*

Pakes et al. [5] alludes to (i) and it is an immediate consequence of remarks preceding Example 6. Part (ii) follows from Implication (iv) because the supports of the components in (1) coincide since they have the same jump distribution; only their rate parameters differ.

We end this section with the following formal extension. Suppose now that the support of $Y_i$ comprises a set of positive integer multiples of $v \in N$. Thus, supp$(Y_i)$ has an arithmetic support with span $v$; that is, $v$ is the greatest common divisor of supp$(Y_i)$. A number theoretic lemma asserts that there exists an integer $j'$ such that $jv \in$ s.g.(supp$(Y_i)$) if $j \geq j'$. This same lemma is used to prove that, in a Markov chain, diagonal $n$-step transition probabilities of an aperiodic state are eventually all positive; see Seneta [6, Lemma A.3]. In this context we say that supp$(X)$ has no gaps if $P\{X = jv\} > 0$ for all $j \in Z_+$. The following result is now clear, part (i) of which appears in Steutel and van Harn [4, page 53] with a different proof and part (ii) is a formal extension of Theorem 5 in Satheesh [2].

TABLE 3

| $N(t) = n = n_1 + n_2$ | 0 | 1 | 2 | 3 | 4 | $\cdots$ |
|---|---|---|---|---|---|---|
| $X(t) = n_1 y_1 + n_2 y_2$ | 0 | 3, 7 | 6, 10, 14 | 9, 21, 13, 17 | 12, 28, 20, 16, 24 | $\cdots$ |

TABLE 4

| $N(t) = n = n_1 + n_2$ | 0 | 1 | 2 | 3 | 4 | $\cdots$ |
|---|---|---|---|---|---|---|
| $X(t) = n_1 y_1 + n_2 y_2$ | 0 | 4, 9 | 8, 13, 18 | 12, 17, 22, 27 | 16, 21, 26, 31, 36 | $\cdots$ |

**Theorem 11.** *Suppose that $X$ is DID with d.f. $F$ and that* supp$(Y_i)$ *is arithmetic with span $\nu$. Then*

(i) supp$(F)$ *has at most finitely many gaps.*

(ii) supp$(F)$ *has no gaps if and if only $F(\{\nu\}) > 0$.*

One point of this extension is that if the distribution of $X$ is gap-free in the sense just defined, then so in an obvious sense is the distribution of $k + X$ for any $k \in N$. With this extended notion of "gap-free," **Example 3** ceases to be a contradiction of Remark 9 in Bose et al. [1]. Note that the term "gap" is not defined in this reference. However, even with our broader notion, taking $X$ as in Examples 7–9, the distribution of $1 + X$ provides counterexamples to Remark 9.

## 3. Some General Considerations

Let $X \geq 0$ be a r.v with *d.f.* $F$ and Laplace-Stieltjes transform $\phi(\theta) = E(e^{-\theta X})$. If $X$ is also ID, then $\phi$ has the canonical form

$$\phi(\theta) = E\left(e^{-\theta X}\right) = \exp\left[-\ell\theta - \int_0^\infty \left(1 - e^{-x\theta}\right) M(dx)\right], \tag{2}$$

where $\ell \geq 0$ is a constant and $M$ is a (Lévy) measure with support supp$(M) \subset [0, \infty)$ and satisfying $\int_0^\infty (1 \wedge x) M(dx) < \infty$.

As above, let $F_n$ denote the $n$'th convolution root of $F$ and $l_n = \sup\{x \geq 0 : F_n(x) = 0\}$ its left extremity, $n \in N$. The following result about left extremities generalizes Lemma 1.1 in Seneta and Vere-Jones [7]. Their result is specific to ID distributions and they prove it using an $\epsilon - \delta$ type of analytic argument. The proof below differs.

**Lemma 12.**

(i) *Let $X \geq 0$ be as above. The left extremity $\ell_X$ of $X$ is given by*

$$e^{-\ell_X} = \lim_{\theta \to \infty} \left[\phi(\theta)\right]^{1/\theta}. \tag{3}$$

(ii) *If $X$ is also ID with Laplace-Stieltjes transform (2), then $\ell_n = \ell/n$, $n \in N$.*

*Proof.*

(i) For $\theta > 0$ the $\theta$-norm of $Z = \exp(-X)$ is $\mathcal{N}_Z(\theta) = (EZ^\theta)^{1/\theta} = (\phi(\theta))^{1/\theta}$. A well-known consequence of Jensen's inequality is that $\mathcal{N}_Z(\theta)$ is nondecreasing in

$\theta$ (named Lyapunov's inequality by some authors). Hence, the limit asserted in (i) exists, and it equals $\mathcal{N}_Z(\infty) = $ a.s. sup $Z = e^{-\ell}$; see Loève [8, page 162].

(ii) Referring to (2) and defining the Lévy tail mass function $\overline{M}(\nu) = M(\nu, \infty)$, we have

$$\theta^{-1} \int_0^\infty \left(1 - e^{-\theta x}\right) M(dx) = \int_0^\infty \int_0^x e^{-\theta \nu} d\nu M(dx)$$
$$= \int_0^\infty e^{-\theta \nu} \overline{M}(\nu) d\nu, \tag{4}$$

and the monotone convergence theorem implies that the right-hand side converges to zero as $\theta \to \infty$. The assertion follows from (2). $\qquad\square$

It follows that $\ell_X = \ell$, thus explaining our notation for the so-called drift parameter in (2). We see also that the $F_n$ have a common support only if $\ell = 0$. We assume this condition in what follows and it is equivalent to the inclusion $0 \in$ supp$(F)$. These equivalent conditions thus generalize **Definition 2**. In addition, it follows from (2) and monotone convergence that $F_n(0) = \exp(-n^{-1}\overline{M}(0+))$. This is zero if and only if $M$ is an infinite measure in which case $\ell_M = 0$; that is, $M$ charges all intervals $(0, \epsilon)$ for each $\epsilon > 0$. The following lemma is a restatement of a known result.

**Lemma 13.** *If $\ell = 0$, then* supp$(F_n) = [0, \infty)$ *if and only if $\ell_M = 0$.*

The direct assertion is proved by Tucker [9] as his Case (i). The "only if" assertion is obvious because if $\ell_M > 0$; then supp$(F_n)$ omits at least the interval $(0, \ell_M)$, a gap. **Lemma 13** follows from the general result asserting that if $\{X(t) : t \geq 0\}$ denotes the Lévy process determined by $M$, that is, $X(1) \stackrel{d}{=} X$, then supp$(X(t)) = \{0\} \cup \overline{s.g.}($supp$(M))$. See Steutel and van Harn [4, pages 110 and 111] for an essentially analytic proof. **Lemma 13** is intuitively obvious from the fact that the sample paths of $\{X(t)\}$ increase by jumps whose sizes and rates are controlled by $M$. So if $\ell_M = 0$, then the Lévy process can hit any subinterval of $(0, \infty)$ with positive probability. This follows in essence because jumps can be arbitrarily small, and hence no interval can be "jumped across" by almost all sample paths.

It is worth noting that if $M$ is an infinite measure and hence that supp$(F_n) = [0, \infty)$, then all $F_n$ are continuous and they can be either absolutely continuous or singular continuous (or have both as components). This follows from small changes to constructions such as in Orey [10].

So gaps are possible in the general case only for compound Poisson distributions. The following example complements Examples 7, 8, and 9 by showing that if $M$ charges an interval, no matter how small, then there are only finitely many gaps, possibly none.

*Example 14.* Choose numbers $0 < \delta < c < \infty$ and suppose that supp$(M) = [c, c + \delta]$. Thus,

$$s.g. \left( \text{supp} \left( M \right) \right) = \cup_{k \geq 1} \left[ ck, (c + \delta) \, k \right]. \tag{5}$$

This set has gaps but only finitely many because, for all sufficiently large integers $k$, we have $\delta k > c$; that is, $(c + \delta)k > c(k+1)$, meaning that the gaps $\mathscr{G}_k = ((c+\delta)k, c(k+1)) = \emptyset$. In addition, the lengths $|\mathscr{G}_k| = c - \delta k$ of nonempty gaps decrease as $k$ grows.

More generally, if supp$(M)$ has positive Lebesgue measure, then $s.g.(\text{supp}(M))$ contains an interval. It thus follows from Example 14 that $s.g.(\text{supp}(M))$ contains a semi-infinite interval and hence so does supp$(F_n)$ for all $n$. Finally, discrete analogues can be constructed by choosing $c$ and $c + \delta$ above to be rational and then supposing that $M$ assigns positive measure to each rational member of $[c, c + \delta]$.

## Conflict of Interests

The authors declare that there is no conflict of interests regarding the publication of this paper.

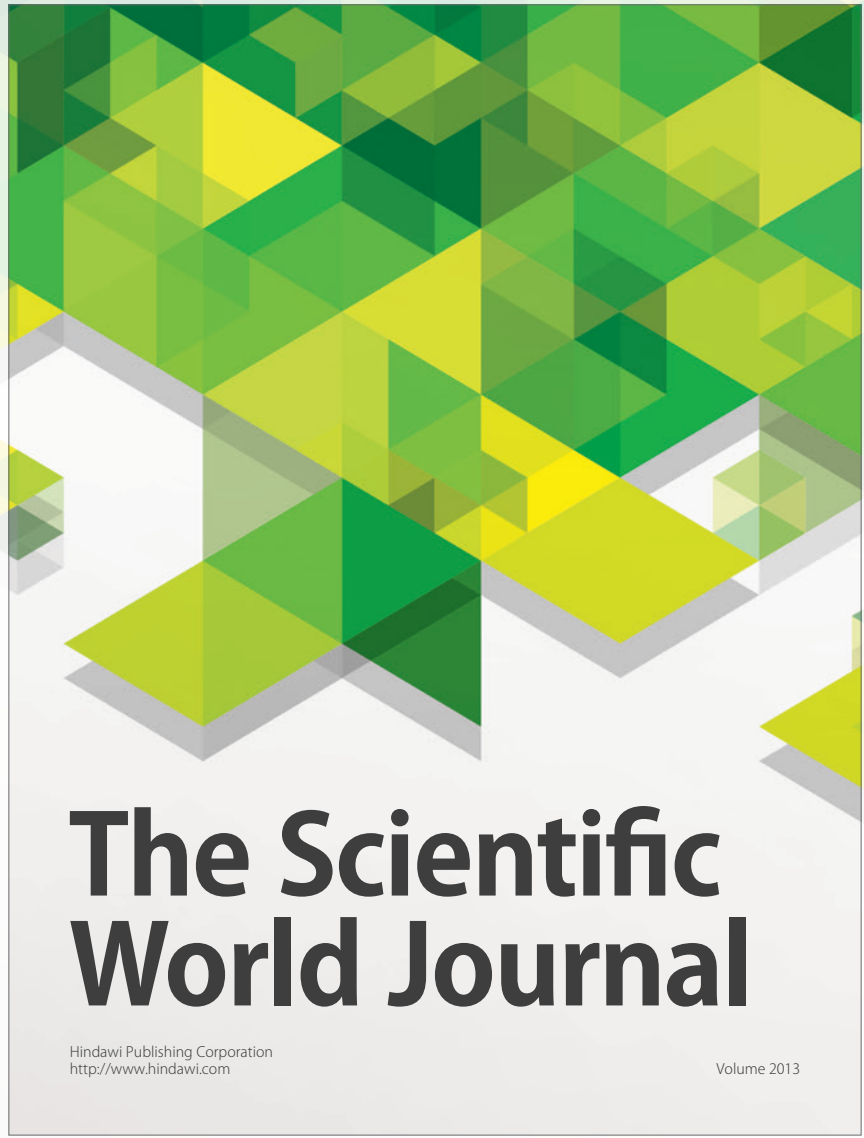

# The Scientific World Journal

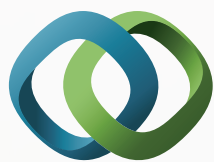

Hindawi

- ▶ Impact Factor **1.730**
- ▶ **28 Days** Fast Track Peer Review
- ▶ All Subject Areas of Science
- ▶ Submit at http://www.tswj.com